\documentstyle[12pt,amssymb,amstex]{amsart}
\numberwithin{equation}{section}

\def\ce{{\cal E}}

\def\co{{\cal O}}

\def\car{{\cal R}}
\def\cs{{\cal S}}

\def\cx{{\cal X}}

\def\ga{{\frak A}}

\def\gd{{\frak D}}
\def\ge{{\frak E}}

\def\gam{{\frak M}}

\def\bc{{\mathbb C}}

\def\bm{{\mathbb M}}
\def\bn{{\mathbb N}}

\def\bq{{\mathbb Q}}
\def\br{{\mathbb R}}

\def\bz{{\mathbb Z}}

\def\a{\alpha}
\def\b{\beta}
  \def\G{\Gamma}
\def\d{\delta}  

\def\e{\epsilon}
\def\l{\lambda} \def\L{\Lambda}

\def\m{\mu}

\def\r{\rho}
\def\s{\sigma} 
\def\t{\tau}
\def\f{\varphi} 
\def\th{\theta}  
\def\om{\omega} \def\Om{\Omega}

\newtheorem{thm}{Theorem}[section]
\newtheorem{lem}[thm]{Lemma}

\newtheorem{defin}[thm]{Definition}

\def\di{\mathop{\rm d}}
\def\sp{\mathop{\rm sp}}
\def\spec{\mathop{\rm spec}}
\def\id{\mathop{\rm id}}
\def\tr{\mathop{\rm Tr}}
\newcommand{\ty}[1]{\mathop{\rm {#1}}}
\newcommand{\nn}{\nonumber}

\begin{document}

\title[diagonalizability of Markov states]
{diagonalizability of non homogeneous quantum Markov states and
associated von Neumann algebras}
\author{Francesco Fidaleo}
\address{Francesco Fidaleo\\
Dipartimento di Matematica\\
Universit\`{a} di Roma ``Tor Vergata''\\
Via della Ricerca Scientifica, 00133 Roma, Italy}
\email{{\tt fidaleo@@mat.uniroma2.it}}
\author{Farruh Mukhamedov}
\address{Farruh Mukhamedov\\
Department of Mechanics and Mathematics\\
National University of Uzbekistan\\
Vuzgorodok, 700095, Tashkent, Uzbekistan} \email{{\tt
far75m@@yandex.ru}}

\begin{abstract}
We clarify the meaning of diagonalizability of quantum Markov
states. Then, we prove that each non homogeneous quantum Markov
state is diagonalizable. Namely, for each Markov state $\f$ on the
spin algebra
$\gam:={\displaystyle\overline{\bigotimes_{j\in\bz}\bm_{d_{j}}(\bc)}^{\,C^{*}}}$
there exists a suitable maximal Abelian subalgebra $\gd\subset\gam$, a
Umegaki conditional expectation $\ge:\gam\mapsto\gd$ and a
Markov measure $\m$ on $\spec(\gd)$ such that
$\f=\f_{\m}\circ\ge$, the Markov state $\f_{\m}$, being the state
on $\gd$ arising from the measure $\m$. An analogous result is
true for non homogeneous quantum processes based on the forward or
the backward chain. Besides, 
we determine the type of the von Neumann factors generated by
GNS representation associated with translation invariant or periodic 
quantum Markov states.
\vskip 0.3cm \noindent
{\bf Mathematics Subject Classification}: 82A15, 46L35, 46L50, 82B20, 
60J99.\\
{\bf Key words}: Mathematical quantum statistical mechanics; 
Classification of von Neumann factors; Non commutative measure, integration and
probability; Lattice systems; Quantum Markov processes.
\end{abstract}

\maketitle

\section{introduction}

It is known that, in quantum statistical mechanics, concrete
systems are identified with states on corresponding algebras. In
many cases, the algebra can be chosen to be a quasi--local algebra
of observables. The states on these algebras satisfying
Kubo--Martin--Schwinger boundary condition, as is known, 
describe equilibrium states of the
quantum system under consideration. On the other hand, 
for classical systems with the
finite radius of interaction, limiting Gibbs measures are know to
be Markov random fields, see e.g. \cite{D, KSK, Pr}. In connection with this, 
there is a
problem to construct analogues of non commutative Markov
chains, which arise from quantum statistical
mechanics and quantum field theory in a natural way. This problem was 
firstly explored in \cite{Ac} by introducing quantum 
Markov chains on the algebra of quasi--local observables. In the last 
decades, the investigation of quantum Markov processes had a 
considerable growth, in view of natural applications to quantum 
statistical mechanics, quantum field theory and quantum information 
theory. The reader is referred to \cite{AFi}--\cite{ALV}, \cite{BKJW}, 
\cite{FNW} and the references cited therein, for recent development of 
the theory of quantum stochastic processes and their applications. 

The investigation of a particular class of quantum Markov chains, 
called quantum Markov states,
was pursued in \cite{AF, ALi}, where connections with 
properties of the modular operator of the states under consideration 
were established. This provides natural applications to temperature 
states arising from suitable quantum spin models, that is natural 
connections with the KMS boundary 
condition.\footnote{Most of the states arising from Markov processes
considered in \cite{FNW} describe ground states 
(i.e. states at zero temperature) of certain models of quantum spin 
chains.}

In \cite{AFi}, the most general
one dimensional quantum Markov state has been considered. Among the 
other results concerning the structure of such states, the connection with
classes of local Hamiltonians satisfying certain commutation
relations and quantum Markov states has been obtained. The situation 
arising from quantum Markov states on the chain, describes one 
dimensional  models of statistical mechanics with mutually commuting 
nearest neighbour interactions. Namely, one dimensional quantum Markov 
states are very near to be (diagonal liftings of) ``Ising type'' models, 
apart from 
noncommuting boundary terms, see Section 6 of \cite{AFi}. 

In the present paper, we clarify
the meaning of diagonalizability of one dimensional non homogeneous
quantum Markov states. Namely, in Section 3 we prove that each non homogeneous
quantum Markov state is diagonalizable, that is, for each Markov state
$\f$ on the spin algebra
$\gam:={\displaystyle\overline{\bigotimes_{j\in\bz}\bm_{d_{j}}(\bc)}^{\,C^{*}}}$
there exists a suitable maximal Abelian subalgebra $\gd\subset\gam$
(called {\it diagonal} in the sequel), a
Umegaki conditional expectation $\ge:\gam\mapsto\gd$ and a
Markov measure $\m$ on $\spec(\gd)$ such that
$\f=\f_{\m}\circ\ge$, the Markov state $\f_{\m}$ being the state
on $\gd$ arising from the measure $\m$. 
This allows us also to clarify a question raised 
in Section 6 of  \cite{AFi},
relative to the r\^ole played by the non commuting boundary terms 
naturally arising from quantum Markov states, see Section 4 below.

The first
diagonalizability result for quantum Markov states is contained in
\cite{Sh} for quantum Markov states generated by a Markov operator. In
\cite{GZ}, the diagonalizability of more general one dimensional
translation invariant quantum Markov states on the forward chain was proved, 
without any 
statement about the Markovianity of the underlying classical measure. The proof in
\cite{GZ} of diagonalizability heavily depend on the commuting square
condition \eqref{iss5} for the increasing sequence of Umegaki conditional 
expectations. The proof of \eqref{iss5}, omitted in \cite{GZ}, 
easily follows by the 
direct ispection of
the structure of local expected subalgebras and potentials, the last  
investigated in detail in \cite{AFi}.

It is known that factors of type $\ty{III}$ naturally arise in quantum field 
theory, statistical physics,
representations of groups, see e.g. \cite{L} and the references cited 
therein. Basically, the systematic 
investigation of the type of factors
generated by the GNS representations of states  
naturally appearing in quantum field theory and in quantum statistical 
mechanics, was an interesting problem since the
pioneering work of Araki and Wyss \cite{AWy}. In
\cite{AWo, P}, a family of representations of uniformly
hyperfinite algebras was constructed. They can be treated as 
free quantum lattice systems. 
In this case, most of the factors corresponding to
these representations are of type $\ty{III}$. However, the product states 
can be viewed as Gibbs
states of Hamiltonian systems in which interactions between
particles of the system are absent. So, it is natural to consider
quantum lattice systems with nontrivial interactions, which lead
us to treat firstly Markov
states, as it was mentioned above. Simple examples of such systems are 
the Ising and  Potts
models. The quantum version of the last ones are diagonal liftings 
(i.e. they are constructed in a trivial way from the corresponding 
classical models, see e.g. Section 5) of classical processes. They   
have been studied in several papers, see e.g. 
\cite{GM, MR1, MR2, S} and the reference cited therein. 

Full analysis relative to 
the type of von Neumann 
algebras arising from general Markov states, or even states associated 
to quantum Markov 
processes on multidimensional lattices, is still an open problem. 

In Section 5, we can partially  
solve this problem for physically relevant one dimensional quantum 
Markov states, that 
is for some examples of translation invariant or periodic states.
Namely, Section 5 of  this paper is devoted to determine the
type of von Neumann factors arising from the GNS representations 
associated to quantum Markov states (for the classification of the 
type $\ty{III}$ factors, see \cite{C}). This is done by using the 
simultaneous diagonalizability of the nearest neighbour terms of the 
interaction associated to quantum Markov states. This classification 
result in the corrected form established in Theorem \ref{noma}) seems to be not 
known even for the Ising model, or for states arising from classical 
Markov chains, the last treated in some detail in Section 5.

Contrary to the  
situations present in literature, the states 
considered here appears as 
nondiagonal liftings of  
classical Markov 
processes which are diagonalizable in a nontrivial way by the result 
proven in Section 3, see Section 4 for a discussion about this
point.\footnote{Other nontrivial quantum liftings of classical Markov chains 
are constructed and studied in \cite{AFi3}. Apart from the standard 
applications to statistical mechanics, possible applications to 
quantum information theory are expected for the last processes.} 
However, it
should be also noted that in \cite{SV}, some properties of 
general diagonal state were studied in 
relation to representations of ``large'' groups of unitaries on 
Hilbert spaces, but
concrete constructions of states were not considered there.

\section{preliminaries}

We start by recalling some well--known facts about inclusions of
finite dimensional $C^{*}$--algebras.

Let $N\subset M$ be an inclusion of
finite dimensional $C^{*}$--algebras. Consider the finite sets
$\{p_{i}\}$, $\{q_{j}\}$ of all the minimal central projections of
$M$, $N$ respectively. We symbolically write
$$
\sum_{j}q_{j}N\subset\sum_{i}p_{i}M\,.
$$

Let us set $M_i:=M_{p_i}$, $N_j:=N_{q_j}$, $M_{ij}:=M_{p_iq_j}$,
$N_{ij}:=N_{p_iq_j}$. Then, we have inclusions $N_{ij}\subset M_{ij}$
of finite dimensional factors. Hence,
\begin{equation}
\label{a}
M_{ij}\sim N_{ij}
\otimes\bar{N}_{ij}
\end{equation}
for other finite dimensional factors
$\bar{N}_{ij}$.\footnote{The square root of the dimension of
$\bar{N}_{ij}$ is precisely the multiplicity of which the piece
$q_{j}N\subset N$ appears into the piece $p_{i}M\subset M$.}

Consider the canonical traces $\tr_{M}$, $\tr_{N}$, that is the
traces which assign unit values on minimal projections. Notice that
$\tr_{M}=\tr_{M}\circ E$ where $E$ is the conditional expectation of
$M$ onto
${\displaystyle \sum_{i,j}q_{j}(p_{i}M)q_{j}}$ given by
$$
E(x)=\sum_{i,j}q_{j}p_{i}xq_{j}\,.
$$

Taking into account
the identification \eqref{a} and the last considerations, one can write
symbolically
\begin{equation*}
\tr{}_{M}=\bigoplus_{i,j}\big(\tr{}_{N_{ij}}\otimes\tr{}_{\bar{N}_{ij}}\big)\,.
\end{equation*}

Further, the completely positive, 
$\big(\tr{}_{M},\tr{}_{N}\big)$--preserving linear map
$E^{M}_{N}$ of $M$ onto $N$ is given by
\begin{equation}
\label{c}
E^{M}_{N}=\bigoplus_{i,j}\big(\id{}_{N_{ij}}\otimes\tr{}_{\bar{N}_{ij}}\big)\,.
\end{equation}

Let $\f$ be a positive functional on $M$, together with its
restriction $\f\lceil_{N}$ to $N$. Consider the corresponding 
Radon--Nikodym
derivatives $T^{\f}_{M}$, $T^{\f}_{N}$ w.r.t. the canonical traces
$\tr_{M}$, $\tr_{N}$ respectively. We get
\begin{equation}
\label{d}
T^{\f}_{N}=E^{M}_{N}(T^{\f}_{M})\,.
\end{equation}

The starting point of our analysis is the $C^{*}$--infinite tensor product
$$
\gam:=\overline{\bigotimes_{j\in\bz}M_{j}}^{\,C^{*}}
$$
where for $j\in\bz$, 
\begin{equation}
\label{aaalg}
M_{j}=\bm_{d_{j}}(\bc)\,.
\end{equation} 

With an abuse of
notations, we denote with the same symbols elements of local
algebras, and their canonical embeddings into bigger (local)
algebras if this cause no confusion. For $k\leq l$, we denote by
$M_{[k,l]}$ the local algebra relative to the segment
$[k,l]\subset\bz$. Let $\cs(\gam)$ be the set of all states
on $\gam$. The restriction of a state $\f\in\cs(\gam)$ to
$M_{[k,l]}$ will be denoted by $\f_{[k,l]}$.

Suppose we have an increasing sequence
$\{N_{[k,l]}\}_{k\leq l}$ of local algebras such that
\begin{align*}
&N_{[k,k]}\subset M_{[k,k]}\equiv M_{k}\,,\quad
N_{[k,k+1]}\subset M_{[k,k+1]}\,,\\
&M_{[k,l]}\subset N_{[k-1,l+1]}\subset M_{[k-1,l+1]}\,,\quad k\leq l\,.
\end{align*}

Consider an increasing
sequence of $C^{*}$--algebras
$\{D_{[k,l]}\}_{k\leq l}$ where $D_{[k,l]}$ is maximal Abelian in
$N_{[k,l]}$. 

A {\it diagonal
algebra} $\gd\subset\gam$ is the Abelian $C^{*}$--subalgebra of $\gam$
obtained as
$$
\gd:=\overline{\big(\lim_{\stackrel{\longrightarrow}
{[k,l]\uparrow\bz}}D_{[k,l]}\big)}^{\,C^{*}}
$$
for $D_{[k,l]}$, $N_{[k,l]}$ as above.

We deal only with {\it locally faithful} states (i.e. states on
$\gam$ with faithful restrictions to local subalgebras), even if most 
of the forthcoming analysis applies to non faithful states as well. For
$\f\in\cs(\gam)$, locally faithful,
the {\it generalized conditional expectation}, or $\f$--expectation,
$\e^{\f}_{k,l}:M_{[k,l+1]}\mapsto M_{[k,l]}$ is the completely
positive $\f$--preserving linear map associated to the inclusion
$M_{[k,l]}\subset M_{[k,l+1]}$
defined in \cite{AC}. We refer the reader to that paper
for the precise definition and further details on the 
Accardi--Cecchini 
generalized conditional expectations.

\section{diagonalizability of Markov states}

Let $\f\in\cs(\gam)$ be a locally faithful state.
\begin{defin}
The state $\f\in\cs(\gam)$ is said to be a Markov state if,
for $k,l\in\bz$, $k<l$, we have
\begin{equation*}
\e^{\f}_{k,l}\lceil_{M_{[k,l-1]}}
=\id{}_{M_{[k,l-1]}}\,.
\end{equation*}
\end{defin}

Quantum Markov states was firstly studied in \cite{Ac, AF}. They are
relevant in quantum  statistical mechanics. The
structure of quantum Markov states was intensively studied in
\cite{AFi, ALi}, where most of their properties were understood.

We briefly report useful results relative to the structure of Markov
states. We refer the reader to \cite{AFi} for details and proofs.

After taking the ergodic limit of the $\f$--expectations
$\e^{\f}_{k,l}$, and a decreasing martingale convergence theorem
(\cite{AFi}, Section 5), it
is possible to recover a sequence $\{\ce^{j}\}_{j\in\bz}$ of
transition expectations which are Umegaki conditional expectations
$$
\ce^{j}:M_{j}\otimes M_{j+1}\mapsto R_{j}\subset M_{j}
$$
such that
\begin{align}
\label{aa}
&\f_{[k,l]}(A_{k}\otimes\cdots\otimes A_{l})\\
=&\f_{[k,k]}(\ce^{k}(A_{k}\otimes\cdots\otimes\ce^{l-1}(A_{l-1}
\otimes A_{l})\cdots))\nn
\end{align}
for every $k,l\in\bz$ with $k<l$, and
$A_{k}\otimes\cdots\otimes A_{l-1}\otimes A_{l}$
any linear generator of $M_{[k,l]}$.
Let $\{P^{j}_{\om_{j}}\}_{\om_{j}\in\Om_{j}}$ be the set of all minimal
central projections of the range $R_{j}:=\car(\ce^{j})$ of
$\ce^{j}$. Put
$$
B_{j}:=\sum_{\om_{j}\in\Om_{j}}P^{j}_{\om_{j}}M_{j}P^{j}_{\om_{j}}\,,
$$
and
$$
B_{[k,l]}:=\bigoplus_{k\leq j\leq l}B_{j}\,.
$$

Consider the conditional expectation $E^{j}:M_{j}\mapsto B_{j}$ given by
$$
E^{j}(A):=\sum_{\om_{j}\in\Om_{j}}P^{j}_{\om_{j}}aP^{j}_{\om_{j}}\,.
$$

Define
\begin{equation}
\label{aaa}
E_{[k,l]}:=\bigoplus_{k\leq j\leq l}E^{j}\,.
\end{equation}

By \eqref{aa}, it is easy to show that
$$
\f_{[k,l]}=\f_{[k,l]}\circ E_{[k,l]}\,.
$$

After the identification
$M_{j}{}_{P^{j}_{\om_{j}}}\cong P^{j}_{\om_{j}}M_{j}P^{j}_{\om_{j}}$
(i.e. the reduced algebra $M_{j}{}_{P^{j}_{\om_{j}}}$
acting on $P^{j}_{\om_{j}}\bc^{d_{j}}$), we have
\begin{equation*}
M_{j}{}_{P^{j}_{\om_{j}}}=N^{j}_{\om_{j}}\otimes\bar N^{j}_{\om_{j}}
\end{equation*}
for finite dimensional factors $N^{j}_{\om_{j}}$, $\bar N^{j}_{\om_{j}}$.
We can write after the last identifications,
\begin{equation}
\label{aab}
B_{[k,l]}:=\bigoplus_{\om_{k},\dots,\om_{l}}
\big(N^{k}_{\om_{k}}\otimes\bar N^{k}_{\om_{k}}\big)
\otimes\cdots\otimes
\big(N^{l}_{\om_{l}}\otimes\bar N^{l}_{\om_{l}}\big)\,.
\end{equation}

Consider the potentials $\{h_{M_{[k,l]}}\}_{k\leq l}$ obtained
by the formula
\begin{equation}
\label{is}
\f_{[k,l]}=\tr{}_{M_{[k,l]}}(e^{-h_{M_{[k,l]}}}\,\cdot\,)\,.
\end{equation}

Then $h_{M_{[k,l]}}$ has the nice decomposition
\begin{equation}
\label{as}
h_{M_{[k,l]}}=\bigoplus_{\om_{k},\dots,\om_{l}}
h^{k}_{\om_{k}}\otimes h^{k}_{\om_{k},\om_{k+1}}
\otimes\cdots\otimes
h^{l-1}_{\om_{l-1},\om_{l}}\otimes\hat h^{l}_{\om_{l}}
\end{equation}
for selfadjoint elements $h^{j}_{\om_{j}}$, $\hat h^{j}_{\om_{j}}$,
$h^{j}_{\om_{j},\om_{j+1}}$ localized in $N^{j}_{\om_{j}}$,
$\bar N^{j}_{\om_{j}}$, $\bar N^{j}_{\om_{j}}\otimes N^{j+1}_{\om_{j+1}}$
respectively.

Defining
\begin{align*}
&H_{j}:=\sum_{\om_{j}}P^{j}_{\om_{j}}
(h^{j}_{\om_{j}}\otimes I)P^{j}_{\om_{j}}\,,\quad
\widehat{H}_{j}:=\sum_{\om_{j}}P^{j}_{\om_{j}}
(I\otimes \hat h^{j}_{\om_{j}})P^{j}_{\om_{j}}\\
&H_{j,j+1}:=\sum_{\om_{j},\om_{j+1}}
(P^{j}_{\om_{j}}\otimes P^{j+1}_{\om_{j+1}})
(I\otimes h^{j}_{\om_{j},\om_{j+1}}\otimes I)
(P^{j}_{\om_{j}}\otimes P^{j+1}_{\om_{j+1}})\,,
\end{align*}
we find sequences of selfadjoint operators $\{H_{j}\}_{j\in\bz}$,
$\{\widehat{H}_{j}\}_{j\in\bz}$ localized in $M_{[j,j]}\equiv
M_{j}$, and $\{H_{j,j+1}\}_{j\in\bz}$ localized in $M_{[j,j+1]}$
respectively, satisfying the commutation relations
\begin{align}
\label{is1}
&[H_{j},H_{j,j+1}]=
[H_{j,j+1},\widehat{H}_{j+1}]\\
=&[H_{j},\widehat{H}_{j}]=[H_{j,j+1},H_{j+1,j+2}]=0\nn\,,
\end{align}
such that
\begin{equation}
\label{is2}
h_{M_{[k,l]}}=H_{k}+\sum_{j=k}^{l-1}H_{j,j+1}+\widehat{H}_{l},
\end{equation}
for each $k\leq l$.

In Section 5 of \cite{AFi} it is proven also the converse. Namely, if
$\f\in\gam$ is locally faithful, with potentials having 
the form \eqref{is2}, for addenda localized as above, and
satisfying the commutation relations \eqref{is1}, then it is a Markov
state.

We are ready to prove the diagonalizability result for quantum Markov
states.
\begin{thm}
\label{mmain} 
Let $\f\in\cs(\gam)$ be a Markov state. Then there
exists a diagonal algebra $\gd\subset\gam$, a classical Markov
process, with Markov measure $\m$  on $\spec(\gd)$ w.r.t. the same
order--localization of $\bz$, and a Umegaki conditional
expectation $\ge:\gam\mapsto\gd$ such that $\f=\f_{\m}\circ\ge$,
where $\f_{\m}$ is the state on $\gd$ corresponding to the measure
$\m$.
\end{thm}
\begin{pf}
Let $R_{j}$ be the range of the (Umegaki) transition expectation
$\ce^{j}$, with relative commutant $R_{j}^{c}:=R_{j}'\bigwedge M_{j}$.
Define
\begin{align*}
&N_{[k,k]}:=Z(R_{k})\,,\quad N_{[k,k+1]}:=R_{k}^{c}\otimes R_{k+1}\,,\\
&N_{[k,l]}:=R_{k}^{c}\otimes M_{[k+1,l-1]}\otimes R_{l}\,,\quad k<l+1\,.
\end{align*}

For each $k\leq j<l$, and $\om_{j}\in\Om_{j}$, choose a maximal
Abelian subalgebra $D^{j}_{\om_{j},\om_{j+1}}$ of
$\bar N^{j}_{\om_{j}}\otimes N^{j+1}_{\om_{j+1}}$ containing
$h^{j}_{\om_{j},\om_{j+1}}$. Put
\begin{align*}
&D_{[k,k]}:=N_{[k,k]}\equiv Z(R_{k})\,,\\
&D_{[k,l]}:=\bigoplus_{\om_{k},\dots,\om_{l}}
\big(D^{k}_{\om_{k},\om_{k+1}}\otimes\cdots\otimes
D^{l-1}_{\om_{l-1},\om_{l}}\big)\,,\quad k<l\,,\\
&\gd:=\overline{\big(\lim_{\stackrel{\longrightarrow}
{[k,l]\uparrow\bz}}D_{[k,l]}\big)}^{\,C^{*}}\,.
\end{align*}

According to our definition, $\gd$ is a diagonal algebra of $\gam$ 
as the $D_{[k,l]}$ are
increasing and
maximal Abelian in the $N_{[k,l]}$. Consider the potentials
$h_{N_{[k,l]}}$ associated to the restrictions $\f\lceil_{N_{[k,l]}}$. We
get by \eqref{d},
$$
e^{-h_{N_{[k,l]}}}=E^{M_{[k,l]}}_{N_{[k,l]}}\big(e^{-h_{M_{[k,l]}}}\big)\,.
$$

Taking into account \eqref{c}, \eqref{as}, we obtain
\begin{equation}
\label{iss2}
h_{N_{[k,l]}}=K_{k}+\sum_{j=k}^{l-1}H_{j,j+1}+\widehat{K}_{l}
\end{equation}
for
\begin{align}
\label{iss3}
&K_{j}:=-\sum_{\om_{j}}\ln\big(\tr{}_{N^{j}_{\om_{j}}}e^{-h^{j}_{\om_{j}}}\big)
P^{j}_{\om_{j}}\,,\\
&\widehat{K}_{j}:=-\sum_{\om_{j}}
\ln\big(\tr{}_{\bar N^{j}_{\om_{j}}}e^{-\hat h^{j}_{\om_{j}}}\big)
P^{j}_{\om_{j}}\nn\,.
\end{align}

Summarizing, by restricting ouselves to the sequence $\{N_{[k,l]}\}_{k\leq l}$,
we find a collection
$\{h_{N_{[k,l]}}\}_{k\leq l}$ of mutually
commuting potentials, with $h_{N_{[k,l]}}\in D_{[k,l]}$, arising from
a nearest neighbour interaction, see \eqref{is1}, \eqref{iss2},
\eqref{iss3}. Namely, $\{h_{N_{[k,l]}}\}_{k\leq l}\subset\gd$.

Let $E_{k,l}:N_{[k,l]}\mapsto D_{[k,l]}$ be the canonical conditional expectation
of $N_{[k,l]}$ onto the maximal abelian subalgebra
$D_{[k,l]}$.\footnote{Let ${\displaystyle M\equiv\sum_{i}p_{i}M}$ be a finite
dimensional
$C^{*}$--algebra, $\{p_{i}\}$ being the set of all its minimal central
projections, and $D\subset M$ a maximal Abelian subalgebra.
Then there exists a complete set of
matrix--units $\{e^{i}_{k_{i},l_{i}}\}$ for $M$ such that $D$ is
generated by the diagonal part $\{e^{i}_{k_{i},k_{i}}\}$. The
canonical expectation $E$ of $M$ onto the diagonal algebra $D$ is
easily given by
$$
E\big(\sum_{i,k_{i},l_{i}}a^{i}_{k_{i},l_{i}}e^{i}_{k_{i},l_{i}}\big)
=\sum_{i,k_{i}}a^{i}_{k_{i},k_{i}}e^{i}_{k_{i},k_{i}}\,.
$$}
We have
\begin{equation}
\label{iss4}
\f\lceil_{N_{[k,l]}}\equiv\tr{}_{N_{[k,l]}}
\big(e^{-h_{N_{[k,l]}}}\,\cdot\,\big)
=\tr{}_{N_{[k,l]}}\big(e^{-h_{N_{[k,l]}}}E_{k,l}(\,\cdot\,)\big)\,.
\end{equation}

Further,
\begin{equation}
\label{iss5}
E_{k-1,l+1}\lceil_{N_{[k,l]}}=E_{k,l}\,.
\end{equation}

Indeed by projectivity,
$$
E_{k,l}=E_{k,l}\circ E_{[k,l]}
$$
with $E_{[k,l]}$ given in \eqref{aaa}. The compatibility condition
\eqref{iss5} immediately
follows taking into account \eqref{aab}.

Let $\f_{\m}:=\f\lceil_{\gd}$, where $\m$ is the probability measure
on $\spec(\gd)$ associated to $\f\lceil_{\gd}$. By \eqref{iss5},
$$
\ge_{0}:=\lim_{\stackrel{\longrightarrow}{[k,l]\uparrow\bz}}E_{k,l}
$$
is well--defined on ${\displaystyle\bigcup_{k,l}N_{[k,l]}}$ (which
is a dense subalgebra of $\gam$), and extends by continuity to a
Umegaki conditional expectation $\ge$ of $\gam$ onto $\gd$.
Further, by \eqref{iss4},
$\f=\f\circ\ge_{0}\equiv\f_{\m}\circ\ge_{0}$ on localized elements
of $\gam$. By a standard continuity argument, we obtain
$\f=\f_{\m}\circ\ge$. The fact that $\m$ is a Markov measure on
$\spec(\gd)$ w.r.t the order--localization of $\bz$ is checked in
the appendix.
\end{pf}

The diagonalizability result for
translation invariant quantum Markov states is contained in
\cite{GZ} for homogeneous processes on the forward chain, without any 
mention about the Markovianity of the underlying classical processes.
As in our situation, the proof of the diagonalizability in Theorem 4.1
of \cite{GZ} heavily depends on the commuting square condition
\eqref{iss5}.\footnote{The proof of \eqref{iss5}, missing in \cite{GZ}, 
would follow by general results contained in Section I.1 of \cite{SV}.}
In the most general situation considered here (hence, including the 
case considered in \cite{GZ}), \eqref{iss5} easily follows by 
a direct ispection of the structure of local
expected subalgebras and potentials investigated in detail in 
\cite{AFi}, and reported in the present paper fpr the convenience of 
the reader.

We end the present section by noticing that an analogous result can be
proven for non homogeneous processes on one--side (forward or backward)
ordered chains. By looking at support--projections of local restritions
of states (or equivalently by defining the Markov property directly in
terms of Umegaki transition expectations, see \cite{AFi}, Definition
2.1), it is straighforward to prove the diagonalizability result for
general (non necessarily locally faithful) Markov states on ordered
chains.

\section{from quantum Markov states to quantum statistical mechanics}

In standard models of statistical mechanics describing classical or
quantum spin systems, one considers on a quasi--local
algebra $\ga$, local
Hamiltonians $\{h_{\L}\}_{\L\subset\bz^{d}}$, $\L$ bounded,
satisfying suitable conditions. Then, one constructs the finite volume
Gibbs states (to simplify matter we reduce ourselves to the case
with inverse temperature $\b=1$)
\begin{equation}
\label{gs}
\f_{\L}:=Z^{-1}\tr{}_{\ga _{\L}}\big(e^{-h_{\L}}\,\cdot\,\big)\,,
\end{equation}
$Z$ being the partition function, see e.g. \cite{BR, Ru, S}. The
local Hamiltonians $h_{\L}$ are usually based on a interaction term
describing the mutual interaction of all spins in the volume $\L$,
and a boundary term arising from some fixed boundary conditions
imposed to spins surronding the bounded region $\L$. After
extending the $\f_{\L}$ to all of $\ga$, each $*$--weak limit
${\displaystyle\lim_{\L_{n}\uparrow\bz^{d}}\f_{\L_{n}}}$ of the
net $\{\f_{\L}\}_{\L\subset\bz^{d}}$ is an infinite volume Gibbs
state, or a DLR state (KMS state in quantum setting) for the system
under consideration see \cite{D2, D3, HHW, LR}.

In the classical case, it is stated for finite range interactions,
that an infinite volume Gibbs state arises from a $\d$--Markov
process and vice--versa, $\d$ being the range of the interaction,
see e.g. \cite{D, KSK, Pr}. For ordered unidimensional chains, a quantum
analogue of that result is proven in \cite{AFi}, provided that the
``leading'' terms $\{H_{j,j+1}\}_{j\in\bz}$ commute with each
other, see also \cite{AFi2} for connected results relative to the 
multidimensional case. In quantum setting, it can happen that
$\{h_{\L}\}_{\L\subset\bz^{d}}$ does not generate a commutative
algebra due to the boudary effects (see \cite{AFi}, Section 6).

In the present paper we have shown that, starting from a quantum Markov state on
${\displaystyle\gam\equiv\overline{\bigotimes_{j\in\bz}M_{d_{j}}}^{\,C^{*}}}$,
we can recover a nontrivial filtration $\{N_{[k,l]}\}_{k\leq l}$
of $\gam$ and an increasing sequence $\{D_{[k,l]}\}_{k\leq l}$ of
Abelian algebras with the $D_{[k,l]}$ non trivial (i.e. not arising 
from the standard tensor structure of $\gam$) maximal Abelian
subalgebras of the $N_{[k,l]}$ such that $\f$ is the lifting
of $\f\lceil_{\gd}$, the last one being a classical
Markov state on 
${\displaystyle\gd:=\overline{\big(\lim_{\stackrel{\longrightarrow}
{[k,l]\uparrow\bz}}D_{[k,l]}\big)}^{\,C^{*}}}$,
constructed by the compatible sequence of U\-me\-ga\-ki conditional expectations
$E_{k,l}:N_{[k,l]}\mapsto D_{[k,l]}$ preserving the canonical trace
$\tr_{N_{[k,l]}}$.\footnote{The restriction of
$\tr_{N_{[k,l]}}$ to $D_{[k,l]}$ is the uniform measure which assignes the same
weight 1 to the minimal projections of $D_{[k,l]}$.}
This is possible as the (nearest neighbour) potentials 
$\{h_{N_{[k,l]}}\}_{k\leq l}$  generate a
commutative subalgebra of $\gd$.

As it is straighforwardly seen, the converse is also true. Namely, one can
start with any fixed filtration $\{N_{[k,l]}\}_{k\leq l}$ as
above, together with a nearest neighbour interaction
\begin{equation}
\label{nne}
h_{k,l}=\sum_{j=k}^{l-1}H_{j,j+1}
\end{equation}
with $\{H_{j,j+1}\}_{j\in\bz}$ mutually commuting. By adding
boundary terms $K_k$ and $\widehat{K}_l$ to \eqref{nne} such that all addenda
commute with each other, one can construct for finite regions 
$\L=[k,l]$, finite volume Gibbs
states $\big\{\f_{\L}\big\}_{\L\subset\bz}$ as in \eqref{gs}, 
associated to the Hamiltonian
$$
h_{[k,l]}=K_k+h_{k,l}+\widehat{K}_l
$$
having the same form as in \eqref{iss2}.

Each $*$--weak limit point of the sequence $\big\{\f_{\L}\big\}_{\L\subset\bz}$,
gives rise to a Markov state on $\gam$
which is the lifting of a classical Markov state on a ``diagonal''
algebra, due to the commutativity of the
$h_{[k,l]}$.\footnote{Notice that, besides the limiting
subsequence $[k_{n},l_{n}]\uparrow\bz$, the thermodinamical limits
might depend also on the chosen boundary terms.}

Taking into account the above considerations, one can assert that each
quantum Markov state on the ordered chain arises from some underlying
(non trivial) classical Ising model.\footnote{For Markov states with
multidimensional indices, where there is no canonical order (i.e. for
the Markov fields considered in \cite{AFi2}), it is expected the appearance of
non diagonalizable examples.}

The quantum character of such states manifests itself in the following
way. In order to construct (or recover) such states, one should take
into account various nontrivial local filtrations of $\gam$, 
together with various (commuting) boundary terms.

If one chooses to investigate quantum Markov states by considering only
the natural filtration $\{M_{[k,l]}\}_{k\leq l}$ of $\gam$, one
obtains a leading term as that in \eqref{nne}.
But non commuting boundary terms could naturally arise
in \eqref{is2}, see the examples in Section 6 of \cite{AFi}.
In the constuctive approach, the appearance of
such non commuting boundary terms cannot be disregarded in order
to obtain general infinite volume Gibbs states for nearest 
neighbour interactions. Here, it should be noted that if the
nearest neighbour model is translation invariant or periodic, 
then according to Theorem 1 of 
\cite{A}, we will have a unique quantum Markov state corresponding to the 
considered model. Namely, for translation invariant or periodic 
models, the construction of quantum Markov states does not
depend on boundary terms.

\section{types of von Neumann algebras associated with quantum Markov 
states}

In this section we investigate the type of von Neumann 
factors generated by the GNS representation associated with the
quantum Markov states.

Let us consider the $C^{*}$--algebra $\gam$ defined in Section 2. 
The shift automorphism of 
the algebra $\gam$ will be
denoted by $\th$. A state  $\f\in \cs(\gam)$ is called 
{\it$l$--periodic} if $\f(\th^l(x))=\f(x)$ for all $x\in\gam$. If
$l=1$, $\f\in \cs(\gam)$ is {\it
translation invariant}. Notice that, in order to have 
$l$--periodicity, it is necessary 
$d_{j+l}=d_j$, $j\in\bz$, for the $d_j$ in \eqref{aaalg}.
Further, we have for localized Hamiltonians \eqref{as}, and their 
leading terms \eqref{nne},
\begin{equation*}
h_{M_{[j+l,k+l]}}=h_{M_{[j,k]}}\,\quad h_{j+l,k+l}=h_{j,k}
\end{equation*}
for all $j,k\in\bz$. In order to avoid the trivial situation, we 
consider only non--tracial locally
faithful translation invariant or $l$--periodic Markov states. 
This means that $h_{0,l}\neq\bc I$, that is $h_{0,l}$ is nontrivial.

We are going to connect the type of the von Neumann factor 
$\pi_{\f}(\gam)''$ with properties of the spectrum 
$\s\big(h_{0,l}\big)$ of the fundamental 
block $h_{0,l}$ of the leading term \eqref{nne} 
of the canonical Hamiltonian 
associated to $\f$. 

Due to commuting properties of the
$h_{M_{[-n,n]}}$ (see \eqref{is1}), the following strong limit
$$
\s^{\f}_t(A)=\lim_{n\to\infty}e^{it
h_{M_{[-n,n]}}}Ae^{-it h_{M_{[-n,n]}}}, \ \ A\in\gam
$$
exists. Further, $\f$ is a KMS (at inverse temperature 1) for 
$\s^{\f}$. According to Theorem 1 of \cite{A}, it is the unique KMS 
state for $\s^{\f}$, and $\pi_{\f}(\gam)''$ is a factor. Notice that we have also
\begin{equation*}
\s^{\f}_t(A)=\lim_{n\to\infty}e^{ith_{-n,n}}Ae^{-ith_{-n,n}}\,.
\end{equation*}

The extension to all of $\pi_{\f}(\gam)''$, denoted 
also by $\s^{\f}$, is precisely the modular group associated to the 
normal extension of $\f$ (denoted also by $\f$) to $\pi_{\f}(\gam)''$. 

Let $\sp(\t)$ be the Arveson spectrum of the action $\t$ of a locally 
compact group on a $C^{*}$--algebra.\footnote{For the 
definition of the Arveson spectrum $\sp(\t)$, as well as 
$\sp{}^{\t}(A)$, see e.g. \cite{Ped}.}
Denote $\s_{t}^{n}:=\text{ad}\big(e^{it h_{-ln,ln}}\big)$, $l$ being the 
period of the state under consideration.
\begin{lem}
\label{ffff}
In the above situation, we have
$$
\sp(\s^{\f})\subset\overline{\bigcup_{n}
\big(\s\big(h_{-ln,ln}\big)-\s\big(h_{-ln,ln}\big)\big)}\,.
$$
\end{lem}
\begin{pf}
By passing to the regrouped algebra, we can consider $l=1$. Taking 
into account the commuting properties of the interaction, we have
\begin{align*}
\sp(\s^{\f})=&\overline{\bigcup_{n}
\bigcup_{A\in M_{[-n,n]}}\sp{}^{\s^{\f}}(A)}
=\overline{\bigcup_{n}
\bigcup_{A\in M_{[-n,n]}}\sp{}^{\s^{n+1}}(A)}\\
\subset&\overline{\bigcup_{n}
\bigcup_{A\in M_{[-n-1,n+1]}}\sp{}^{\s^{n+1}}(A)}
=\overline{\bigcup_{n}
\sp(\s^{n+1}\lceil_{M_{[-n-1,n+1]}})}\,.
\end{align*}

The proof follows by Proposition 14.13 of \cite{St}.
\end{pf}
\begin{lem}
\label{fffff}
 Let $\{x_1,\dots,x_n\}\subset\br\backslash\{0\}$ such that
$x_i/x_j\in\bq$ for all $i,j$. Then 
\begin{equation}
\label{rat}
\{x_1,\dots,x_n\}\subset\bz\ln\a
\end{equation}
for some $\a\in(0,1)$.
\end{lem}
\begin{pf} From our assumptions, we have
$$
\frac{x_1}{x_i}=\frac{p_i}{q_i},  \ \ \ \ i=2,...,n\,,
$$
where $p_i\in\bn\backslash\{0\}$, 
$q_i\in\bz\backslash\{0\}$.\footnote{The best $\a$ in \eqref{rat} is 
the minimum of the $\a\in(0,1)$ such that \eqref{rat} is true. It can be 
obtained by changing the reference element, and compute  
all the corresponding $\a$ in \eqref{rat1} by taking relatively prime pairs 
$p_i,q_i$. The minimum we are looking for, is precisely the smallest 
among all these $\a$.} Define
\begin{equation}
\label{rat1}
\a:=e^{-\frac{|x_1|}{\prod_{j=2}^np_j}}\,.
\end{equation}

Then
\begin{align*}
&x_1=-\text{sign}(x_1)\big(\prod_{j=2}^np_j\big)\ln\a\,,\\
&x_i=-q_i\big(\prod_{\stackrel{j=2}{j\neq i}}^np_j\big)\ln\a\,,\quad 
i=2,\dots,n\,.
\end{align*}
\end{pf}

Let $h_{0,l}$ be the fundamental 
block of the leading term of the canonical Hamiltonian 
associated to the locally faithful Markov state $\f$. 
Consider, for $h,k,h',k'\in\s(h_{0,l})$ with $h\neq k$,
$h'\neq k'$, the following fractions 
$\frac{h-k}{h'-k'}$.
\begin{thm}
\label{noma}
Let $\f\in\cs(\gam)$ be a locally faithful Markov state. The following 
assertions hold true.
\begin{itemize}
\item[(i)] If $\big\{\frac{h-k}{h'-k'}\big\}\subset\bq$, then 
$\pi_{\f}(\gam)''$ is a type $\ty{III_{\l}}$ factor for some 
$\l\in(0,1)$.
\item[(ii)] If $\pi_{\f}(\gam)''$ is a type $\ty{III_{1}}$ factor,
then $\big\{\frac{h-k}{h'-k'}\big\}\nsubseteq\bq$.
\end{itemize}
\end{thm}
\begin{pf}
As before, we can consider only translation invariant Markov 
states. By applying Theorem 3.1 of \cite{Se}, we get for the Connes 
invariant $\G$ (see \cite{C}), 
$\G(\pi_{\f}(\gam)'')\equiv\G(\s^{\f})=\sp(\s^{\f})$. Further, this 
means also that $\pi_{\f}(\gam)''$ is a type $\ty{III_{\l}}$ factor, 
$\l\in(0,1]$, as we are considering non--tracial states. Then, it is 
enough to prove the former, the latter being a direct consequence of 
the former.

Let $\big\{\frac{h-k}{h'-k'}\big\}\subset\bq$ be satisfied. By Lemma 
\ref{fffff}, 
$$
\{h-k\,|\,h,k\in\s(H_{0,1})\}\subset\bz\ln\a
$$ 
for some $\a\in(0,1)$.

From the simultaneous diagonalizability of the $H_{i,i+1}$, we find that
$$
\s(h_{-n,n})\subset\bigg\{\sum_{i=-n}^{n-1}h_i\, \bigg|\, 
h_i\in\s(H_{0,1})\bigg\}.
$$

Then we have
$$
\s(h_{-n,n})-\s(h_{-n,n})\subset\bigg\{\sum_{i=-n}^{n-1}(h_i-k_i)\, \bigg|\,
h_i,k_i\in\s(H_{0,1})\bigg\}\subset\bz\ln\a\,.
$$

According to Lemma \ref{ffff},
we infer that $\sp(\s^{\f})\subset\bz\ln\a$, that is
$\sp(\s^{\f})$ is discrete. Hence, there is a
number $m\in\bn\backslash\{0\}$ such that $\sp(\s^{\f})=\bz\ln\l$, with
$\l:=\a^{m}$.
Thus $\pi_{\f}(\gam)''$ is a type $\ty{III_{\l}}$ factor.
\end{pf} 

Here, it should be noted that one might argue that the spectrum 
$\s(h_{0,l})$ of the fundamental block of the Hamiltonian associated 
to the periodic Markov state $\f$, completely
determines the type of $\pi_{\f}(\gam)''$. Unfortunately, we are not 
able to prove the reverse statements in Theorem \ref{noma}, 
contrarily to that is asserted in literature.\footnote{For example, 
the proof of the connected results in \cite{MR2} seems to be 
incomplete.}

Even if one can costruct by results in Section 4 of \cite{AFi}, a wide 
class of quantum Markov states 
to which the previous results apply, in order to explain some 
natural applications of Theorem \ref{noma} to pre--assigned models,
we are going to consider some natural examples.

\subsection{Ising model}

In this situation,
$$
\gam=\overline{\bigotimes_{\bz}\bm_{2}(\bc)}^{\,C^{*}}\,.
$$

The Ising model on
$\bz$ is defined by the following formal Hamiltonian
$$
H=-\sum_{j\in\bz}J_{j,j+1}\s_z^j\s_z^{j+1},
$$
where $J_{j,j+1}\in\br$ are coupling constants and $\s_z^j$ is the
Pauli matrix $\s_z$ on the $j$--th site. Further, we suppose that
the coupling constants are defined by
$$
J_{j,j+1}= \left\{ \begin{array}{ll} J_1, \ \ \ \textrm{if} \ \
j\in 2\bz,\\[2mm]
J_2, \ \ \ \textrm{if} \ \
j\in 2\bz+1,\\
\end{array}
\right.
$$
where $J_1,J_2\in\br$. It is known (see \cite{A}) that for the
given Hamiltonian there exists a unique Gibbs state $\f$ on 
$\gam$ which is $2$--periodic. 
In this case,
the operators $H_{j,j+1}$ have the following form
$$
H_{j,j+1}=\left\{ \begin{array}{ll}
            \left( \begin{array}{cccc}
             J_1 & 0 & 0 & 0 \\
             0  & -J_1 & 0 & 0\\
             0  & 0 & -J_1 & 0\\
             0  & 0 & 0 & J_1\\
             \end{array}\right), \ \ \ \textrm{if} \ \ j\in 2\bz,\\[12mm]
             \left( \begin{array}{cccc}
             J_2 & 0 & 0 & 0 \\
             0  & -J_2 & 0 & 0\\
             0  & 0 & -J_2 & 0\\
             0  & 0 & 0 & J_2\\
             \end{array}\right), \ \ \ \textrm{if} \ \ j\in
             2\bz+1,\\[3mm]
            \end{array}
        \right.
$$

The
spectrum of $H_{j,j+1}$ is $\{J_1,-J_1\}$ if $j\in 2\bz$,
$\{J_2,-J_2\}$ if $j\in 2\bz+1$ respectively. Now if $J_1/J_2$ is
rational, the rationality condition of Theorem \ref{noma} is 
satisfied, this
infers that the von Neumann factor $\pi_{\f}(\gam)''$ is of type 
$\ty{III_\l}$, for some $\l\in(0,1)$.

\subsection{Markov process} 

Consider a discrete Markov process with the
state space ${\bf d}:=\{1,\dots,d\}$ and the transition probabilities
defined by means of a stochastic matrix
$P=(p_{ij})_{i,j=1}^d$ with (not all equal) $p_{ij}>0$  for all $i,j$. 
Consider the canonical inclusion
$$
\gd=\overline{\bigotimes_{\bz}\bc_{d}}^{\,C^{*}}
\subset\gam=\overline{\bigotimes_{\bz}\bm_{d}(\bc)}^{\,C^{*}}\,.
$$

Here, $\gd\sim C(\Om)$, where $\Om=\prod\limits_{\bz}{\bf d}$. Let 
$\m_{P}$ be the translation invariant 
Markov measure on $\Om$ determined by the transition 
matrix $P$. Define the diagonal lifting of the classical process 
associated to $P$ as
$$
\f(A):=\int_{\Om}\ge(A)(\om)\m_{P}(\di\om)\,,
$$
where $\ge$ is the canonical Umegaki conditional 
expectation of $\gam$ onto the Abelian algebra $\gd$. It is clear
that such state is a translation invariant quantum Markov state.

It is not hard to check that the corresponding $H_{j,j+1}$ operator
has the form
$$
H_{j,j+1}= \left (\begin{array}{cccccc}
B^{(1)}&0&\cdot&\cdot&\cdot&0 \\ 0&B^{(2)}&0&\cdot&\cdot&0\\
\cdot&\cdot&\cdot&\cdot&\cdot&\cdot\\
\cdot&\cdot&\cdot&\cdot&\cdot&\cdot\\
\cdot&\cdot&\cdot&\cdot&\cdot&\cdot\\
0&0&\cdot&\cdot&\cdot&B^{(d)}\\ \end{array} \right), 
$$ 
where
$B^{(k)}=(b_{ij,k})_{i,j=1}^{d}$, $k=1,\dots,d$ are  d$\times$d
diagonal matrices such that 
$$ 
b_{ij,k}= \left\{
\begin{array}{ll} -\ln p_{k,i},\,\,   i=j,\ i=1,\dots,d\\ 
0,\, \qquad \quad i\neq j \end{array} \right.
$$

If there exist integers $m_{ij}$
$i,j\in\{1,\dots,d\}$, and some number $\a\in(0,1)$ such
that $ \frac{p_{11}}{p_{i,j}}=\a^{m_{ij}}$,
then we easily see that the rationality condition of 
Theorem \ref{noma} is satisfied, this means that
the von Neumann factor $\pi_{\f}(\gam)''$ is of type 
$\ty{III_\l}$, for some $\l\in(0,1)$.
This result extends a result of \cite{GM}.

\section{appendix}

For the convenience of the reader, we verify that the measure 
$\m$ on $\gd$ associated to 
$\f\lceil_{\gd}$ is a Markov measure on $\spec(\gd)$ w.r.t. the
order--localization of $\bz$.

For our pourpose,
it suffices to verify that for
every $k\leq n\leq l$ in $\bz$ and $A\in\spec(D_{[k,n]})$ and
$B\in\spec(D_{[n,l]})$
we have for the conditional probability,
$$
P\big(A\cap B\big|\bar\om_{n}\big)=P\big(A\big|\bar\om_{n}\big)
P\big(B\big|\bar\om_{n}\big)\,,
$$
where $\bar\om_{n}$ is a fixed point in $\spec(Z(R_{n}))\equiv\Om_{n}$.
In order
to make computations, we should see the past $D_{[k,n]}$, the
present algebra $D_{[n,n]}\equiv Z(R_{n})$, and the future algebra
$D_{[n,l]}$ inside the ambient algebra $D_{[k,l]}$. In such a
situation
$$
\spec(D_{[k,l]})
=\stackrel{\bf{.}}{\bigcup}_{\om_{k},\dots,\om_{l}}S^{k}_{\om_{k},\om_{k+1}}
\times\cdots\times S^{k}_{\om_{l-1},\om_{l}}\,.
$$

Taking into account Formulae \eqref{as} and \eqref{iss2}, we compute
for
$$
f:=\sum_{\om_{k},\dots,\om_{l}}
\cx_{S^{k}_{\om_{k},\om_{k+1}}\times\cdots\times S^{l-1}_{\om_{l-1},\om_{l}}}
f^{k}_{\om_{k},\om_{k+1}}\otimes\cdots\otimes
f^{l-1}_{\om_{l-1},\om_{l}}\,,
$$
\begin{align}
\label{integ}
&\f(f)=\sum_{\om_{k},\dots,\om_{l}}\bigg(\int_{S^{k}_{\om_{k},\om_{k+1}}}
T^{k}_{\om_{k},\om_{k+1}}f^{k}_{\om_{k},\om_{k+1}}\bigg)\nn\\
&\times\cdots\times\bigg(\int_{S^{l-1}_{\om_{l-1},\om_{l}}}
T^{l-1}_{\om_{l-1},\om_{l}}f^{l-1}_{\om_{l-1},\om_{l}}\bigg)
\end{align}
where the densities $T$ are positive functions, and $\int$ assignes
weight $1$ to atoms.

We start by noticing that, inside $D_{[k,l]}$, we get for
$P^{n}_{\bar\om_{n}}$,
\begin{align}
\label{integ1}
P^{n}_{\bar\om_{n}}=&
\sum_{\hskip-5pt{\renewcommand\arraystretch{.7}\begin{array}c
\scriptstyle \om_{k},\dots,\om_{n-1},\\
\scriptstyle \om_{n+1},\dots,\om_{l}\end{array}}}
\cx_{S^{k}_{\om_{k},\om_{k+1}}}\otimes\cdots\otimes
\cx_{S^{n-1}_{\om_{n-1},\bar\om_{n}}}\nn\\
\otimes&\cx_{S^{n}_{\bar\om_{n},\om_{n+1}}}\otimes\cdots\otimes
\cx_{S^{l-1}_{\om_{l-1},\om_{l}}}\,.
\end{align}

Now, if $A, B\in\spec(D_{[k,l]})$ are localized in the past and in
the future of $n$ respectively, we have inside $D_{[k,l]}$,
\begin{align}
\label{prdt}
\cx_{A}P^{n}_{\bar\om_{n}}=
\sum_{\hskip-5pt{\renewcommand\arraystretch{.7}\begin{array}c
\scriptstyle a\in A|\r(a)=\bar\om_{n}\\
\scriptstyle \om_{n+1},\dots,\om_{l}\end{array}}}
&\cx_{\{a^{k}_{\om_{k}(a),\om_{k+1}(a)}\}}\otimes\cdots\otimes
\cx_{\{a^{n-1}_{\om_{n-1}(a),\bar\om_{n}}\}}\nn\\
\otimes&\cx_{S^{n}_{\bar\om_{n},\om_{n+1}}}\otimes\cdots\otimes
\cx_{S^{l-1}_{\om_{l-1},\om_{l}}}\,,\nn\\
P^{n}_{\bar\om_{n}}\cx_{B}=
\sum_{\hskip-5pt{\renewcommand\arraystretch{.7}\begin{array}c
\scriptstyle \om_{k},\dots,\om_{n-1}\\
\scriptstyle b\in B|\l(b)=\bar\om_{n} \end{array}}}
&\cx_{S^{k}_{\om_{k},\om_{k+1}}}\otimes\cdots\otimes
\cx_{S^{n-1}_{\om_{n-1},\bar\om_{n}}}\\
\otimes&\cx_{\{b^{n}_{\bar\om_{n},\om_{n+1}(b)}\}}\otimes\cdots\otimes
\cx_{\{b^{l-1}_{\om_{l-1}(b),\om_{l}(b)}\}}\,,\nn\\
\cx_{A}P^{n}_{\bar\om_{n}}\cx_{B}=
\sum_{\hskip-5pt{\renewcommand\arraystretch{.7}\begin{array}c
\scriptstyle a\in A|\r(a)=\bar\om_{n}\\
\scriptstyle b\in B|\l(b)=\bar\om_{n} \end{array}}}
&\cx_{\{a^{k}_{\om_{k}(a),\om_{k+1}(a)}\}}\otimes\cdots\otimes
\cx_{\{a^{n-1}_{\om_{n-1}(a),\bar\om_{n}}\}}\nn\\
\otimes&\cx_{\{b^{n}_{\bar\om_{n},\om_{n+1}(b)}\}}\otimes\cdots\otimes
\cx_{\{b^{l-1}_{\om_{l-1}(b),\om_{l}(b)}\}}\nn\,.
\end{align}

Here,
\begin{align*}
&a=a^{k}_{\om_{k}(a),\om_{k+1}(a)}\times\cdots\times
a^{n-1}_{\om_{n-1}(a),\om_{n}(a)}\,,\\
&b=b^{n}_{\om_{n}(b),\om_{n+1}(b)}\times\cdots\times
b^{l-1}_{\om_{l-1}(b),\om_{l}(b)}
\end{align*}
are generic points of $A$, $B$ respectively, and with these notations,
$\r(a):=\om_{n}(a)$, $\l(b):=\om_{n}(b)$.

Taking into account \eqref{integ}, \eqref{integ1} and \eqref{prdt},
we have the following computations.
\begin{align*}
\f\big(P^{n}_{\bar\om_{n}}\big)
=&\sum_{\hskip-5pt{\renewcommand\arraystretch{.7}\begin{array}c
\scriptstyle \om_{k},\dots,\om_{n-1},\\
\scriptstyle \om_{n+1},\dots,\om_{l}\end{array}}}
\big(\int_{S^{k}_{\om_{k},\om_{k+1}}}
T^{k}_{\om_{k},\om_{k+1}}\big)
\times\cdots\times\big(\int_{S^{n-1}_{\om_{n-1},\bar\om_{n}}}
T^{n-1}_{\om_{n-1},\bar\om_{n}}\big)\\
\times&
\big(\int_{S^{n}_{\bar\om_{n},\om_{n+1}}}
T^{n}_{\bar\om_{n},\om_{n+1}}\big)
\times\cdots\times\big(\int_{S^{l-1}_{\om_{l-1},\om_{l}}}
T^{l-1}_{\om_{l-1},\om_{l}}\big)\\
\equiv&\bigg(\sum_{\om_{k},\dots,\om_{n-1}}
\big(\int_{S^{k}_{\om_{k},\om_{k+1}}}
T^{k}_{\om_{k},\om_{k+1}}\big)
\times\cdots\times\big(\int_{S^{n-1}_{\om_{n-1},\bar\om_{n}}}
T^{n-1}_{\om_{n-1},\bar\om_{n}}\big)\bigg)\\
\times&\bigg(\sum_{\om_{n+1},\dots,\om_{l}}
\big(\int_{S^{n}_{\bar\om_{n},\om_{n+1}}}
T^{n}_{\bar\om_{n},\om_{n+1}}\big)
\times\cdots\times\big(\int_{S^{l-1}_{\om_{l-1},\om_{l}}}
T^{l-1}_{\om_{l-1},\om_{l}}\big)\bigg)\,.
\end{align*}

\begin{align*}
&\f\big(\cx_{A}P^{n}_{\bar\om_{n}}\cx_{B}\big)\\
=\sum_{\hskip-5pt{\renewcommand\arraystretch{.7}\begin{array}c
\scriptstyle a\in A|\r(a)=\bar\om_{n}\\
\scriptstyle b\in B|\l(b)=\bar\om_{n} \end{array}}}
&T^{k}_{\om_{k}(a),\om_{k+1}(a)}\big(a^{k}_{\om_{k}(a),\om_{k+1}(a)}\big)
\times\cdots\times
T^{n-1}_{\om_{n-1}(a),\bar\om_{n}}
\big(a^{n-1}_{\om_{n-1}(a),\bar\om_{n}}\big)\\
\times&T^{n}_{\bar\om_{n},\om_{n+1}(b)}
\big(b^{n}_{\bar\om_{n},\om_{n+1}(b)}\big)\times\cdots\times
T^{l-1}_{\om_{l-1}(b),\om_{l}(b)}
\big(b^{l-1}_{\om_{l-1}(b),\om_{l}(b)}\big)\\
\equiv\bigg(\sum_{a\in A|\r(a)=\bar\om_{n}}
&T^{k}_{\om_{k}(a),\om_{k+1}(a)}\big(a^{k}_{\om_{k}(a),\om_{k+1}(a)}\big)
\times\cdots\times
T^{n-1}_{\om_{n-1}(a),\bar\om_{n}}
\big(a^{n-1}_{\om_{n-1}(a),\bar\om_{n}}\big)\bigg)\\
\times\bigg(\sum_{b\in B|\l(b)=\bar\om_{n}}
&T^{n}_{\bar\om_{n},\om_{n+1}(b)}
\big(b^{n}_{\bar\om_{n},\om_{n+1}(b)}\big)\times\cdots\times
T^{l-1}_{\om_{l-1}(b),\om_{l}(b)}
\big(b^{l-1}_{\om_{l-1}(b),\om_{l}(b)}\big)\bigg)\,.
\end{align*}

\begin{align*}
&\f\big(\cx_{A}P^{n}_{\bar\om_{n}}\big)\\
=\sum_{\hskip-5pt{\renewcommand\arraystretch{.7}\begin{array}c
\scriptstyle a\in A|\r(a)=\bar\om_{n}\\
\scriptstyle \om_{n+1},\dots,\om_{l}\end{array}}}
&T^{k}_{\om_{k}(a),\om_{k+1}(a)}\big(a^{k}_{\om_{k}(a),\om_{k+1}(a)}\big)
\times\cdots\times
T^{n-1}_{\om_{n-1}(a),\bar\om_{n}}
\big(a^{n-1}_{\om_{n-1}(a),\bar\om_{n}}\big)\\
\times
&\big(\int_{S^{n}_{\bar\om_{n},\om_{n+1}}}
T^{n}_{\bar\om_{n},\om_{n+1}}\big)
\times\cdots\times\big(\int_{S^{l-1}_{\om_{l-1},\om_{l}}}
T^{l-1}_{\om_{l-1},\om_{l}}\big)\\
\equiv\bigg(\sum_{a\in A|\r(a)=\bar\om_{n}}
&T^{k}_{\om_{k}(a),\om_{k+1}(a)}\big(a^{k}_{\om_{k}(a),\om_{k+1}(a)}\big)
\times\cdots\times
T^{n-1}_{\om_{n-1}(a),\bar\om_{n}}
\big(a^{n-1}_{\om_{n-1}(a),\bar\om_{n}}\big)\bigg)\\
\times\bigg(\sum_{\om_{n+1},\dots,\om_{l}}
&\big(\int_{S^{n}_{\bar\om_{n},\om_{n+1}}}
T^{n}_{\bar\om_{n},\om_{n+1}}\big)
\times\cdots\times\big(\int_{S^{l-1}_{\om_{l-1},\om_{l}}}
T^{l-1}_{\om_{l-1},\om_{l}}\big)\bigg)\,.
\end{align*}

\begin{align*}
&\f\big(P^{n}_{\bar\om_{n}}\cx_{B}\big)\\
=\sum_{\hskip-5pt{\renewcommand\arraystretch{.7}\begin{array}c
\scriptstyle \om_{k},\dots,\om_{n-1}\\
\scriptstyle b\in B|\l(b)=\bar\om_{n} \end{array}}}
&\big(\int_{S^{k}_{\om_{k},\om_{k+1}}}
T^{k}_{\om_{k},\om_{k+1}}\big)
\times\cdots\times\big(\int_{S^{n-1}_{\om_{n-1},\bar\om_{n}}}
T^{n-1}_{\om_{n-1},\bar\om_{n}}\big)\\
\times&T^{n}_{\bar\om_{n},\om_{n+1}(b)}
\big(b^{n}_{\bar\om_{n},\om_{n+1}(b)}\big)\times\cdots\times
T^{l-1}_{\om_{l-1}(b),\om_{l}(b)}
\big(b^{l-1}_{\om_{l-1}(b),\om_{l}(b)}\big)\\
\equiv\bigg(\sum_{\om_{k},\dots,\om_{n-1}}
&\big(\int_{S^{k}_{\om_{k},\om_{k+1}}}
T^{k}_{\om_{k},\om_{k+1}}\big)
\times\cdots\times\big(\int_{S^{n-1}_{\om_{n-1},\bar\om_{n}}}
T^{n-1}_{\om_{n-1},\bar\om_{n}}\big)\bigg)\\
\times\bigg(\sum_{b\in B|\l(b)=\bar\om_{n}}
&T^{n}_{\bar\om_{n},\om_{n+1}(b)}
\big(b^{n}_{\bar\om_{n},\om_{n+1}(b)}\big)\times\cdots\times
T^{l-1}_{\om_{l-1}(b),\om_{l}(b)}
\big(b^{l-1}_{\om_{l-1}(b),\om_{l}(b)}\big)\bigg)\,.
\end{align*}

Collecting together the last computations, we get
\begin{align*}
P\big(A\cap&B\big|\bar\om_{n}\big)
\equiv\frac{\f\big(\cx_{A}P^{n}_{\bar\om_{n}}\cx_{B}\big)}
{\f\big(P^{n}_{\bar\om_{n}}\big)}\\
=&\frac{\f\big(\cx_{A}P^{n}_{\bar\om_{n}}\big)}
{\f\big(P^{n}_{\bar\om_{n}}\big)}
\frac{\f\big(P^{n}_{\bar\om_{n}}\cx_{B}\big)}
{\f\big(P^{n}_{\bar\om_{n}}\big)}\\
\equiv&P\big(A\big|\bar\om_{n}\big)
P\big(B\big|\bar\om_{n}\big)
\end{align*}
which is the assertion.

\section*{acknowledgements}

We thank L. Accardi for fruitful discussions.
The second--named author (F. M.) acknowledges the Italian CNR for providing
financial support, and Universit\`a di Roma "Tor Vergata" for all
facilities, and in particular L. Accardi and E. Presutti for kind hospitality.


\begin{thebibliography}{9999}

\bibitem{Ac} Accardi L, {\it On noncommutative Markov property},
Funct. Anal. Appl. {\bf 8} (1975), 1--8.

\bibitem{AC} Accardi L., Cecchini C.
{\it Conditional expectations in von Neumann algebras and a theorem
of Takesaki}, J. Funct. Anal. {\bf 45} (1982), 245--273.

\bibitem{AFi} Accardi L., Fidaleo F.
{\it Non homogeneous quantum Markov states and quantum Markov fields},
J. Funct. Anal. {\bf 200} (2003), 324--347.

\bibitem{AFi2} Accardi L., Fidaleo F.
{\it Quantum Markov fields},
Infin. Dimens. Anal. Quantum Probab. Relat. Top. {\bf 6} (2003),
123--138.

\bibitem{AFi3} Accardi L., Fidaleo F.
{\it Entangled Markov chains},
Ann. Mat. Pura Appl., to appear.

\bibitem{AF} Accardi L., Frigerio A.
{\it Markovian cocycles}, Proc. R. Ir. Acad. {\bf 83} (1983), 251--263.

\bibitem{ALi} Accardi L., Liebscher V.
{\it Markovian KMS states for one dimensional spin chains},
Infin. Dimens. Anal. Quantum Probab. Relat. Top. {\bf 2} (1999), 645--661.

\bibitem{ALV} Accardi L., Lu Y. G., Volovich I. {\it Quantum theory 
and its stochastic limit}, Springer, Berlin--Heidelberg--New
York, 2002.

\bibitem{A} Araki H. {\it  On uniqueness of KMS states of one--dimensional quantum
lattice systems}, Commun. Math. Phys. {\bf 44} (1975), 1--7.

\bibitem{AWo} Araki H., Woods E. J. {\it A classification of
factors}, Publ. RIMS Kyoto Univ. {\bf 3} (1968), 51--130.

\bibitem{AWy} Araki H., Wyss W. {\it  Representations of canonical
anticommutations relations}, Helv. Phys. Acta {\bf 37} (1964),
136--159.

\bibitem{BKJW}
Bratteli O., Jorgensen P. E. T., Kishimoto A., Werner R. F. 
{\it Pure states on $\co_d$}, J. Operator Theory {\bf 43} (2000), 97--143. 

\bibitem{BR} Bratteli O., Robinson D. W. {\it Operator algebras and
quantum statistical mechanics, Voll.  I, II}, Springer, Berlin--Heidelberg--New
York, 1981.

\bibitem{C} Connes A., {\it Une classification des facteurs de type III},
Ann. Scient. \'Ec. Norm. Sup. {\bf 6} (1973), 133--252.

\bibitem{D} Dobrushin R. L. {\it The description of a random field by
means of conditional probabilities and conditions of its regularity},
Theor. Probab. Appl. {\bf 13} (1968), 197--224.

\bibitem{D2} Dobrushin R. L. {\it Gibbsian random fields for lattice
systems with pairwaise interactions},
Funct. Anal. Appl. {\bf 2} (1968), 292--301.

\bibitem{D3} Dobrushin R. L. {\it The problem of uniqueness of a
Gibbsian random field and the problem of phase transitions},
Funct. Anal. Appl. {\bf 2} (1968), 302--312.

\bibitem{FNW} Fannes M., Nachtergaele B., Werner R. F. {\it Finitely 
correlated states of quantum spin chains}, Commun. Math. Phys. 
{\bf 144} (1992), 443--490.

\bibitem{GM} Ganikhodjaev N. N., Mukhamedov F. M. {\it  Markov states
on quantum lattice systems and its representations}, Methods 
Funct. Anal. Top. {\bf 4} (1998), 33--38.

\bibitem{GZ} Golodets Ya V., Zholtkevich G. N.
{\it Markovian KMS states},
Theor. Math. Phys. {\bf 56} (1983), 686--690.

\bibitem{HHW} Haag R., Hugenholtz N. M., Winnink M.
{\it On the equilibrium states in quntum statistical mechanics},
Commun. Math. Phys. {\bf 5} (1967), 215--236.

\bibitem{KSK} Kemeny J. G., Snell J.L., Knapp A. W. 
{\it Denumerable Markov chains}, 
Springer--Verlag, Berlin Heidelberg New York 1976.

\bibitem{LR} Lanford III O. E., Ruelle D. {\it Observables at infinity
and states with short range correlations in statistical mechanics},
Commun. Math. Phys. {\bf 13} (1969), 194--215.

\bibitem{L} Longo R. {\it Algebraic and modular structure of von
Neumann algebras of Physics}, Proc. Symp. Pure Math. {\bf 38}
(1982), 551--566.

\bibitem{MR1} Mukhamedov F. M., Rozikov U. A. {\it  von Neumann algebra
corresponding to one phase of inhomogeneous Potts model on a
Cayley tree},  Theor. Math. Phys. {\bf 126} (2001), 169--174.

\bibitem{MR2} Mukhamedov F. M., Rozikov U. A. {\it On Gibbs measures of models
with competing ternary and binary interactions and corresponding
von Neumann algebras}, J. Stat. Phys. {\bf 114} (2004), 
825--848.

\bibitem{Ped} Pedersen G. K.
{\it $C^{*}$--algebras and their automorphism groups}, Academic
press, London, 1979.

\bibitem{P} Powers R. {\it  Representation of uniformly hyperfinite algebras
and their associated von Neumann rings}, Ann. Math. {\bf 81} (1967), 138--171.

\bibitem{Pr} Preston C. J. {\it Gibbs states on countable sets}
Cambridge University Press, London, 1974.

\bibitem{Ru} Ruelle D., {\it Statistical mechanics},
Benjamin, Amsterdam--New York, 1969.

\bibitem{Sh} Shukhov A. G. {\it Entropy of diagonal states on
von Neumann algebras}, Funct. Anal. Appl. {\bf 14} (1980),
163--165.

\bibitem{S}  Sinai Ya G. {\it Theory of phase transitions: Rigorous Results}
Pergamon press, Oxford, 1982.

\bibitem{Se} St{\o}rmer E.
{\it Spectra of states, and asymptotically Abelian
$C^{*}$--algebras}, Commun. Math. Phys. {\bf 28} (1972), 279--294.

\bibitem{St} Str\v{a}til\v{a} S. {\it Modular theory in operator algebras}, 
Abacus press, Tunbridge Wells, Kent, 1981.

\bibitem{SV} Str\v{a}til\v{a} S., Voiculescu D. {\it Representations of
$AF$--algebras and of the group $U(\infty)$}, Springer, Berlin--Heidelberg--New
York, 1975.

\end{thebibliography}
\end{document}